\documentclass{article}
\usepackage[dvips]{color}
\usepackage{amssymb}
\usepackage{amsmath}
\usepackage{hyperref}

\def\F{\mathbb F}

\def\C{\mathbb C}
\def\R{\mathbb R}
\def\S{\mathbb S}
\def\Q{\mathbb Q}
\def\A{\mathbb A}

\def\O{\mathbb O}
\ifx\pdfoutput\undefined
\usepackage{graphicx}
\else
\usepackage[pdftex]{graphicx}
\fi

\def\0{\underline{0}}

\usepackage{tikz}
\usepackage{amsfonts}

\newtheorem{theorem}{Theorem}[section]
\newtheorem{lemma}{Lemma}[section]

\newtheorem{example}{Example}[section]
\newtheorem{remark}{Remark}[section]
\newtheorem{definition}{Definition}[section]

\begin{document}

\begin{center}
{\fontsize{14}{20}\bf
Counting Extended Irreducible Goppa Codes
}
\end{center}

\begin{center}
\textbf{Kondwani Magamba$^{1,2}$ and John A. Ryan$^3$}\\[.25cm]
$^1$Malawi University of Science and Technology, Malawi\\ 
$^{2}$Mzuzu University, Malawi\\
$^{3}$Chombe Boole Research Center,Malawi

\end{center}

%\maketitle

\begin{abstract}
We obtain an upper bound on the number of extended irreducible q-ary Goppa codes of degree $r$ and length $q^n+1$, where $q=p^t$ and $n$ and $r>2$ are prime numbers. 
\end{abstract}


%\IEEEpeerreviewmaketitle

\section{Introduction}

%\IEEEPARstart{T}{he} problem of counting non-extended irreducible Goppa codes has been considered by a number of authors. In \cite{Ryan2}, an upper bound on the number of irreducible Goppa codes is given. Recently, there have been attempts to enumerate extended irreducible Goppa codes. For example, \cite{Ryan4} gives an upper bound on the number of extended irreducible binary quartic Goppa codes of length $2^n+1$. In this paper we obtain an upper bound on the number of extended irreducible $q$-ary Goppa codes of degree $r$ and length $q^n+1$, where $q=p^t$ and $n$ and $r>2$ are prime numbers. To do this, we first develop tools for counting extended irreducible Goppa codes.  

This paper studies extended irreducible Goppa codes. Our interest in Goppa codes stems from the fact they are of practical value as they form the backbone of the McEliece cryptosystem. In the McEliece cryptosystem, one chooses a random Goppa code as a key hence it is important that we know the number of Goppa codes for any given set of parameters. This will help in the assessment of how secure the McEliece cryptosystem is against an enumerative attack. An enumerative attack on the McEliece cryptosystem finds all Goppa codes for the given set of parameters and tests their equivalence with the public code. Research has clearly shown that many Goppa codes become equivalent when extended by a parity check (see \cite{Musukwa} and \cite{Ryan4}) and so it has been suggested that an enumerative attack can be mounted through extended Goppa codes. In this paper we obtain an upper bound on the number of extended irreducible $q$-ary Goppa codes of degree $r$ and length $q^n+1$, where $q=p^t$ and $n$ and $r>2$ are prime numbers. Methods similar to the ones used in the present paper have been used in \cite{Musukwa} and \cite{Ryan4} to count the number of extended irreducible Goppa codes of length $2^n+1$ of degree $2^m$ and degree 4, respectively. However, in this work, we obtain our count by exploiting the action of the projective linear group on the set of elements of degree $r$.

\section{Preliminaries}
We begin this section by defining an irreducible Goppa code.
\subsection{Irreducible Goppa codes}
\begin{definition}
Let $n$ be a positive integer, $q$ be a power of a prime number and $g(z)\in \F_{q^n}[z]$ be irreducible of degree $r$. Let
$L=\F_{q^n}=\{\zeta_i: 0\leq i \leq q^n-1 \}$. Then an irreducible Goppa code $\Gamma(L,g)$ is defined as the set of all vectors $\underline{c}=(c_0,c_1,\ldots,c_{q^n-1})$ with components in $\F_q$ which satisfy the condition
\begin{equation}
\sum_{i=0}^{q^n-1}\frac{c_i}{z-\zeta_i}\equiv 0~ \mbox{mod}~g(z).
 \end{equation}
\end{definition}
The polynomial $g(z)$ is called the Goppa polynomial. Since $g(z)$ is irreducible and of degree $r$ over $\F_{q^n}$, $g(z)$ does not have any root in $L$ and the code is called an irreducible Goppa code of degree $r$. In this paper $g(z)$ is
always irreducible of degree $r$ over $\F_{q^n}$.\\

It can be shown, see \cite{Chen}, that if $\alpha$ is any root of the Goppa polynomial $g(z)$ then  $\Gamma(L, g)$ is completely described by any root $\alpha$ of $g(z)$ and a parity check matrix $\bf{H}(\alpha)$ is given by
\begin{equation}
\bf{H}(\alpha)=\left(\frac{1}{\alpha-\zeta_{0}}\frac{1}{\alpha-\zeta_{1}}\cdots \frac{1}{\alpha-\zeta_{q^n-1}}\right)
\end{equation}
where $L=\F_{q^n}=\{\zeta_i: 0\leq i \leq q^n-1 \}$. This code is usually denoted $C(\alpha)$. 
\subsection{Extended Irreducible Goppa codes} 
Next we give the definition of an extended irreducible Goppa code.
\begin{definition}
Let $\Gamma(L,g)$ be an irreducible Goppa code of length $q^{n}$. Then the extended code $\overline{\Gamma(L,g)}$ is defined by $\overline{\Gamma(L,g)}=\{(c_{0},c_{1},...,c_{q^{n}}): (c_{0},c_{1},...,c_{q^{n}-1})\in \Gamma(L,g)\hspace{0.2cm} \mbox{and}~ \sum_{i=0}^{q^{n}}c_{i}=0\}$.
\end{definition}

\subsection{Enumeration of Matrices of a Given Order in $GL(2,q^n)$}\label{Matriculation}
In this section we obtain an enumeration of matrices of a given order $k$ in $GL(2,q^n)$, where $GL(2,q^n)$ is the group of $2\times 2$-invertible matrices over $\F_{q^n}$. Invertible matrices of a given order will be a powerful tool that will be used in the enumeration of extended irreducible Goppa codes. We first review concepts in matrix algebra. 

We know that a $2\times 2$ matrix $A$ is conjugate (or similar) to a $2\times 2$ matrix $B$ if there exists a non-singular $2\times 2$ matrix $P$ such that $A=P^{-1}BP$. It is well known that the conjugacy of square matrices defines an equivalence relation on and partitions $GL(2,q^n)$ into disjoint sets of equivalence classes called conjugacy classes. We also know that conjugate matrices have the same order.              

Thus the number of matrices in $GL(2,q^n)$ of order $k$ is given by the product of the number of conjugacy classes containing elements of order $k$ and the number of elements in each conjugacy class. Now to find the number of conjugacy classes in $GL(2,q^n)$ we use the fact that all matrices in a given conjugacy class have the same minimal polynomial. Since conjugacy classes partition $GL(2,q^n)$ it follows that the number of distinct minimal polynomials will give us the number of conjugacy classes. We explain this below.

\begin{definition}
Let $A\in GL(2,q^n)$. Then the characteristic polynomial of a matrix $A$ is defined by $\chi_A(x)=det(A-xI_2)$ where $I_2=\left(
\begin{array}{cc}
1&0\\
0&1
\end{array}\right)$ and $det(B)$ is the determinant function. 
\end{definition}
%\begin{example}
%Let $A=\left(
%\begin{array}{cc}
%1&0\\
%0&1
%\end{array}\right)\in GL(2,3)$
%\end{example}
\begin{definition}
We say that $A\in GL(2,q^n)$ is a matrix root of the polynomial\\ $f(x)=a_kx^k+\cdots+a_0\in \F_{q^n}[x]$ if $$f(A)=a_kA^k+\cdots+a_0I_2=0,$$ where $A^i$ is the ith power of $A$ under matrix multiplication.  
\end{definition} 

\begin{definition}
The minimal polynomial of $A\in GL(2,q^n)$ denoted $m_A(x)$ is the monic polynomial $m_A(x)\in \F_{q^n}[x]$ of least degree having $A$ as a matrix root. 
\end{definition}
By the Cayley-Hamilton Theorem, $m_A(x)$ is a factor of $\chi_A(x)$. Since elements in a conjugacy class have the same order and minimal polynomial, minimal polynomials are in a one-to-one correspondence with conjugacy classes in $GL(2,q^n)$. Now, the matrices $A$ are of size $2\times 2$ so $\chi_A(x)$ is a quadratic polynomial, as such $m_A(x)$ is either linear, or a product of two like or unlike linear factors or an irreducible monic quadratic polynomial, see \cite{Alperin}. In this work we are interested in matrices whose minimal polynomial is quadratic and irreducible over $\F_{q^n}$. The following lemma gives the relationship between such minimal polynomials and conjugacy classes of elements of $GL(2,q^n)$. 
\begin{lemma}\label{lemma_T4}
Let $A\in GL(2,q^n)$. If $m_A(x)=x^2-\xi x-\zeta$, where $\xi,\zeta \in \F_{q^n}$, is irreducible then $A$ is conjugate with $\left(
\begin{array}{cc}
0&1\\
\zeta&\xi
\end{array}\right)$ and there are $\frac{q^n(q^n-1)}{2}$ conjugacy classes of length $q^n(q^n-1)$ each.

{\it Proof.} See \cite{Alperin} and \cite{Basheer}.
\end{lemma} 
%We use Lemma \ref{lemma_T4} when we consider the action of $G$ on $\O$.
 
Now to find the number of matrices of order $k$ in $GL(2,q^n)$ we exploit the fact that if $A\in GL(2,q^n)$ has finite order $k$, then its minimal polynomial $m_A(x)$ divides $x^k-1 \in \F_{q^n}[x]$, see \cite{Koo}. So to find all possible minimal polynomials we consider the factorization of $x^k-1$ over $\F_{q^n}$ which can be found in \cite{Lidl} and is stated below.

\begin{theorem}[Lidl and Niederreiter, 1983]{\cite{Lidl}}\label{Lidl}
Let $k$ be a positive integer and $\F_{q^n}=\F_{p^{nt}}$. If $p\nmid k$, then $$x^k-1=\prod_{s\mid k} Q_s(x),$$ where $Q_s(x)$ is the $k-\mbox{th}$ cyclotomic polynomial over $\F_{q^n}$. Furthermore, $Q_k(x)$ factors into $\phi(k)/d$ distinct monic irreducible polynomials of the same degree $d$, where $d$ is the least positive integer such that $q^{nd}\equiv 1 \pmod{k}$.
\end{theorem}
Now, if $k$ is prime then $x^k-1=(x-1)Q_k(x)$. It is easy to see that if there are irreducible quadratic factors in the factorization of $x^k-1$ then they will occur in the factorization of $Q_k(x)$.  

In the following theorem we find the number of matrices $A$ of order $k$ in $GL(2,q^n)$ whose minimal polynomial $m_A(x)$ is an irreducible quadratic polynomial over $\F_{q^n}$.
\begin{theorem}\label{Thm_Yanga}
Let $k$ be a positive integer such that $(q^n,k)=1$ and $k \mid (q^{n}+1)$ but $k \nmid (q^n-1)$. Then the number of matrices $A$ of order $k$ in $GL(2,q^n)$, where $m_A(x)$ is an irreducible quadratic polynomial over $\F_{q^n}$, is $\frac{\phi(k)q^n(q^n-1)}{2}$.
\\ {\rm {\it Proof.}~%We exploit the fact that if $A\in GL(2,q^n)$ has finite order $k$, then the minimal polynomial of $A$, $m_A(x)$ divides $x^k-1$, see \cite{Koo}. The matrices we are interested in are of the form $A=\left(
%\begin{array}{cc}
%a&b\\
%c&d
%\end{array}\right)\in GL(2,q^n)$ where $m_A(x)$ is irreducible over $\F_{q^n}$ so we count irreducible quadratic polynomials in the factorisation of  $x^k-1\in \F_{q^n}[x]$.
 Let $\rho$ be the number of irreducible quadratic polynomials in the factorization of the cyclotomic polynomial $Q_k(x)$ where $(q^n,k)=1$ and $k \mid (q^{n}+1)$ but $k \nmid (q^n-1)$. By Theorem \ref{Lidl}, the cyclotomic polynomial $Q_k(x)$ factors into $\rho=\frac{\phi(k)}{2}$ distinct irreducible quadratic factors. Since each minimal polynomial gives rise to a conjugacy class, the number of conjugacy classes containing elements of order $k$ is equal to $\rho$. From Lemma \ref{lemma_T4}, we know that there are $q^n(q^n-1)$ matrices in a conjugacy class in this case. Thus the number of matrices in $GL(2,q^n)$ of order $k$ whose minimal polynomial $m_A(x)$ is an irreducible quadratic polynomial over $\F_{q^n}$  is $\frac{\phi(k)q^n(q^n-1)}{2}$.}
\end{theorem}
 
\begin{example}\label{Ex_1}
Let's take $q=3$, $n=3$ and $r=7$. Suppose we want to find the number of matrices $A$ of order $7$ in $GL(2,3^3)$, where $m_A(x)$ is an irreducible quadratic polynomial over $\F_{3^3}$. Here $k=7$ and $7|(3^3+1)$. Clearly, $\rho=3$ and $\mu=3^3\times (3^3-1)=702$. So there are $3\times 702=2,106$ matrices in $GL(2,3^3)$ of the required form of order $7$.    
\end{example}

\subsection{Tools for Counting Extended Irreducible Goppa codes}
\subsubsection{The set $\S=\S(n,r)$}
An irreducible Goppa code can be defined by any root of its Goppa polynomial. As such the set of all roots of such polynomials is important and we make the following definition.

\begin{definition}
The set $\S=\S(n,r)$ is the set of all elements in $\F_{q^{nr}}$ of degree $r$ over $\F_{q^{n}} $.  
\end{definition}

\subsubsection{Maps on $\S$}\label{maps_on_S}
We define the following maps on $\S$. The action of the groups arising from these maps will help us to count the number of irreducible Goppa codes and their extended versions.
\begin{enumerate}
\item $\sigma^i:\alpha \to \alpha^{q^i}$ where $\sigma$ denotes the Frobenius
automorphism of $\F_{q^{nr}}$ leaving $\F_q$ fixed and $0 \leq i<nr$.\label{uyu1}
\item $\pi_A:\alpha\to a\alpha+b$ where $A=\left(
\begin{array}{cc}
a&b\\
0&1
\end{array}\right)
\in GL(2,q^n)$.\label{uyu2}
\item $\pi_{B}:\alpha \to \frac{a\alpha+b}{c\alpha+d}$, where $B=\left(
\begin{array}{cc}
a&b\\
c&d
\end{array}\right)
\in GL(2,q^n)$.\label{uyu3} 
\end{enumerate}
\begin{remark}
The composition of maps \ref{uyu1} and \ref{uyu2} sends irreducible Goppa codes into equivalent irreducible Goppa codes and \ref{uyu1} and \ref{uyu3} sends extended irreducible Goppa codes into equivalent extended irreducible Goppa codes (see \cite{Berger}).
\end{remark}

\subsubsection{Groups Acting on $\S$}
%In this section we show that $GL(2,q^n)$ acts on $\S$. We also define the following groups which act on $\S$. The action of these groups on $\S$ will help us count the number of irreducible Goppa codes and their extended versions.
\begin{definition}
Let $G$ denote the set of all maps $\{\sigma^{i}:1\leq i\leq nr\}$. $G$ is a group under the composition of mappings. It is the group of Frobenius automorphisms. It is shown in \cite{Ryan3} that $G$ acts on $\S$.
\end{definition}

\begin{definition}
Let $F$ denote the set of all maps $\left\lbrace\pi_{A}:A=\left(
\begin{array}{cc}
a&b\\
0&1
\end{array}\right)
\in GL(2,q^n) \right\rbrace$. $F$ is a group under the composition of mappings. It is the affine group of linear transformations.
\end{definition}

\begin{definition}
We define the set of all maps $\left\lbrace\pi_{B}:B=\left(
\begin{array}{cc}
a&b\\
c&d
\end{array}\right)
\in GL(2,q^n)\right\rbrace$. This set of maps is a group under the composition of mappings and is isomorphic to the projective linear group which is denoted by $PGL(2,q^n)$.
\end{definition}
Next we show that $PGL(2,q^n)$ acts on $\S$. 
Denote by $[A]$ the image of $A=\left(
\begin{array}{cc}
a&b\\
c&d
\end{array}\right)
\in GL(2,q^n)$ in $PGL(2,q^n)$. Then, for $\alpha \in \S$ and $[A]\in PGL(2,q^n)$ we define the map $[A](\alpha)=\frac{a\alpha+b}{c\alpha+d}$.
Denote by $I_2$ the $2 \times 2$ identity matrix $\left(
\begin{array}{cc}
1&0\\
0&1
\end{array}\right)
$. Then we have $[I_2](\alpha)=\frac{1\alpha+0}{0\alpha+1}=\alpha$.\\
Also, suppose $B=\left(
\begin{array}{cc}
a_1&b_1\\
a_2&b_2
\end{array}\right)$ and $C=\left(
\begin{array}{cc}
a_3&b_3\\
a_4&b_4
\end{array}\right)$ then
 
 \begin{align*}
 [B]([C](\alpha))=&[B]\left(\frac{a_3\alpha+b_3}{a_4\alpha+b_4}\right)\\
 =&\frac{a_1\left(\frac{a_3\alpha+b_3}{a_4\alpha+b_4}\right)+b_1}{a_2\left(\frac{a_3\alpha+b_3}{a_4\alpha+b_4}\right)+b_2}\\
 =&\frac{(a_1a_3+a_4b_1)\alpha+a_1b_3+b_1b_4}{(a_2a_3+a_4b_2)\alpha+a_2b_3+b_2b_4}\\
 =&[BC](\alpha).
 \end{align*}

Thus, the map $[A](\alpha)=\frac{a\alpha+b}{c\alpha+d}$ defines a group action of $PGL(2,q^n)$ on $\S$. This can also be found in \cite{Stich}.

\subsubsection{Actions of $F$, $PGL(2,q^n)$ and $G$}\label{Actions}
We first consider the action of the affine group $F$ on $\S$. For each $\alpha \in \S$, the action of $F$ on $\S$ induces orbits denoted $A(\alpha)$ where $A(\alpha)=\{a\alpha+b:a \neq 0, b \in \F_{q^n}\}$, and called the affine set containing $\alpha$. We denote the set of all affine sets, $ \{A(\alpha):\alpha \in \S \}$, by $\A$. Since $|A(\alpha)|=q^n(q^n-1)$ then $|\A|=|\S|/q^n(q^n-1)$. It can be shown that $G$ acts on the set $\A$, see \cite{Ryan2}. We will then consider the action of $G$ on $\A$ to obtain orbits in $\S$ of $FG$. The number of orbits in $\S$ under $FG$ will then give us an upper bound on the number of irreducible Goppa codes.

Next we consider the action of the group $E=PGL(2,q^n)$ on $\S$. The action of $PGL(2,q^n)$ on $\S$ induces orbits denoted by $O(\alpha)$ where $O(\alpha)=\{\frac{a\alpha +b}{c\alpha +d}: a,b,c,d\in \mathbb{F}_{q^{n}}, ad-bc \neq 0\}$. We will refer to $O(\alpha)$ as a projective linear set. Next we calculate the cardinality of $O(\alpha)$.
\begin{theorem}
For any $\alpha \in \S$, $|O(\alpha)|=q^{3n}-q^n$.\label{Card_O}\\
{\it Proof.} We know from the orbit-stabilizer theorem that  $|O(\alpha)|=|E|/|E_\alpha|$ where $E_\alpha$ is the stabilizer of $\alpha$ in $E=PGL(2,q^n)$. Observe that $E_\alpha$ is trivial, for if some $[B]\in E$, where $B=\left(
\begin{array}{cc}
a&b\\
c&d
\end{array}\right)
$, fixes $\alpha \in \S$ then $[B](\alpha)=\alpha$. That is, $\frac{a\alpha+b}{c\alpha+d}=\alpha$. So, $a\alpha+b=c\alpha^2+d\alpha$, this implies $c\alpha^2+(d-a)\alpha-b=0$. Since the minimal polynomial of $\alpha$ over $\F_{q^n}$ is of degree $r \geq 3$, we conclude that $b=c=0$ and $a=d$, hence $[B]=[I_2]$. As such, $|O(\alpha)|=|E|=q^{3n}-q^n$.
 \end{theorem}

We denote the set of all projective linear sets in $\mathbb{S}$ under the action of $PGL(2,q^n)$ by {\rm $\mathbb{O}$. That is, $\mathbb{O}=\{O(\alpha): \alpha \in \mathbb{S}\}$. Observe that $\mathbb{O}$ partitions the set $\mathbb{S}$} and that $G$ acts on the set $\mathbb{O}$ \cite{Ryan3}.

It is shown in \cite{Ryan3} that each projective linear set $O(\alpha)$ in $\mathbb{O}$ can be partitioned into $q^{n}+1$ affine sets. See the theorem below.

\begin{theorem}\label{O}
For $\alpha \in \mathbb{S}, O(\alpha)=A(\alpha)\cup A(\frac{1}{\alpha})\cup A(\frac{1}{\alpha + 1})\cup A(\frac{1}{\alpha + \xi _{1}}) \cup A(\frac{1}{\alpha + \xi _{2}})\cup \dots \cup A(\frac{1}{\alpha + \xi _{q^{n}-2}})$ where $\mathbb{F}_{q^{n}}=\{0,1,\xi_{1}, \xi_{2},\ldots, \xi_{q^{n}-2}\}$.
\end{theorem}

Observe that the sets $\mathbb{O}$ and $\mathbb{A}$ are different. $\O$ and $\mathbb{A}$ are both partitions of $\mathbb{S}$ but $|\mathbb{A}|=(q^{n}+1) \times |\mathbb{O}|$. 

%$1-\ell$

We will use the actions of $PGL(2,q^n)$ and $G$ on $\S$ to find an upper bound on the number of extended irreducible Goppa codes. Firstly, we will apply the action of $PGL(2,q^n)$ on $\S$ to obtain projective linear sets $O(\alpha)$. Then we will consider the action of $G$ on $\O$. The number of orbits in $\O$ under the action of $G$ will give an upper bound on the  number of extended irreducible Goppa codes. To find the number of orbits we will use the Cauchy-Frobenius theorem which is stated below.
\begin{theorem}\label{Cauchy-Frobenius}
Let $G$ be a finite group acting on a set $X$. For any $g\in G$, let $X(g)$ denote the set of elements of $X$ fixed by $g$. Then the number of orbits in $X$ under the action of $G$ is $\frac{1}{|G|}\sum_{g\in G}|X(g)|$.
\end{theorem}

\subsubsection{Factorization of the polynomial $F_s(x)=cx^{q^s+1}+dx^{q^s}-ax-b$}\label{F_s}
Another tool we shall use when it comes to counting extended irreducible Goppa codes is counting how many roots of the polynomial $F_s(x)=cx^{q^s+1}+dx^{q^s}-ax-b\in \F_{q^n}[x]$ where $ad-bc \neq 0$ lie in $\S$. We will do this by counting the number of irreducible polynomials of degree $r$ in the factorization of $F_s(x)=cx^{q^s+1}+dx^{q^s}-ax-b$. This problem was considered in \cite{Stich} for $F_s(x)=cx^{q^s+1}+dx^{q^s}-ax-b\in \F_{q}[x]$, where $ad-bc \neq 0$ and $s\geq 1$. 
\begin{theorem}[Stichtenoth, H., and Topuzo$\breve{\mbox{g}}$lu, 2012]{\cite{Stich}}\label{Stich_Bho}
Let $A=\left(
\begin{array}{cc}
a&b\\
c&d
\end{array}\right)\in GL(2,q)$ which is not a multiple of the identity matrix. Let the order of $[A]$ in $PGL(2,q)$ be $D$. Then the irreducible factors of $F_s(x)=cx^{q^s+1}+dx^{q^s}-ax-b\in \F_{q}[x]$, $s\geq 1$ are as follows:
\begin{enumerate}
\item irreducible factors of degree $Ds$,
\item irreducible factors of degree $Dk$ with $k<s$, $s=km$ and $\mbox{gcd}(m,D)=1$,
\item irreducible factors of degree $\leq 2$.
\end{enumerate}

\end{theorem} 

\section{Extended Irreducible Goppa codes of degree $r$ and length $q^n+1$, where $n=r$ is prime}
In this section we obtain an upper bound on the number of extended irreducible Goppa codes of degree $r$ and length $q^n+1$, where $n=r$ is prime. We first obtain $\S(n,n)$.
\subsection{$\S(n,r)$ where $n=r$ is prime} 
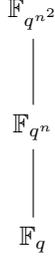
\begin{figure}[htb]
\centering
\begin{tikzpicture}[node distance=1.5cm] 
\node (F1) {$\F_q$}; 
\node (Fn) [above of=F1]{$\F_{q^{n}}$};
\node (Fnr) [above of=Fn]{$\F_{q^{n^2}}$};

\draw (F1) -- (Fn);
\draw (Fn) -- (Fnr);
\end{tikzpicture}
\caption{Lattice of subfields of $\F_{q^{n^2}}$}
\label{242}
\end{figure}

We use a lattice of subfields, as proposed in \cite{Magamba}, to show where elements of $\S(n,n)$ lie and find $|\S(n,n)|$. Figure~\ref{242} shows a lattice of subfields corresponding to $q$ and $n=r$.

\begin{remark}
The number of elements of degree $n$ over $\F_{q^n}$ is $|\S(n,n)|=q^{n^2}-q^{n}.$ 
\end{remark}

\subsection{Action of $G$ on $\A$}
In this section we find the number of affine sets fixed by the subgroups $\langle\sigma\rangle$, $\langle\sigma^{n}\rangle$ and $\langle\sigma^{n^2}\rangle$ of $G$. It is easy to see that the trivial subgroup $\langle\sigma^{n^2}\rangle$, containing the identity, fixes every affine set in $\A$. By an argument similar to the one in  Section III of \cite{Ryan4} corresponding to the actions of the subgroups  $\langle\sigma^n\rangle$ and $\langle\sigma\rangle$ we obtain the following result:
\begin{theorem}\label{Theorem_Puliya_1}
{\rm The number of affine sets fixed under $\langle\sigma^{n}\rangle$ and $\langle\sigma\rangle$ is
$$\left\{ 
\begin{array}{r l}
n-1,& \quad \mbox{if}~ n\neq p~\mbox{and}~n\mid q^n-1 \\
0,& \quad \mbox{if}~ ~ n\neq p~\mbox{and}~n\nmid q^n-1\\
1,& \quad \mbox{if}~ n=p
\end{array} \right..$$
}
\end{theorem}

\subsection{Action of G on $\O$}
In this section we obtain an upper bound on the number of extended irreducible Goppa codes of degree $r$ and length $q^n+1$ where  $n=r$ is prime. As already stated, we will do this  by considering the action of $G$ on $\O$ and applying the Cauchy-Frobenius Theorem. We begin by finding the number of projective linear sets $O(\alpha)$ which are in $\O$. By Theorem \ref{Card_O}, $|O(\alpha)|=q^{3n}-q^n$. We also know that $|\S(n,n)|=q^{n^2}-q^{n}$. Therefore, there are $\frac{q^{n^2}-q^{n}}{q^{3n}-q^n}$ projective linear sets  in $\O$. The number of orbits in $\O$ under the action of $G$ will give us an upper bound on the number of extended irreducible Goppa codes. Clearly, the trivial subgroup $\langle\sigma^{n^2}\rangle$, containing the identity, fixes every projective linear set in $\O$. We now consider the actions of $\langle\sigma^{n}\rangle$ and $\langle\sigma\rangle$. 
\subsubsection{Action of $\langle\sigma^{n}\rangle$, a subgroup of order $n$}\label{Action_O_n_n=r}
Suppose $O(\alpha)\in \O$ is fixed under $\langle\sigma^{n}\rangle$. Then $\langle\sigma^{n}\rangle$ acts on $O(\alpha)=A(\alpha)\cup A(\frac{1}{\alpha})\cup A(\frac{1}{\alpha + 1})\cup A(\frac{1}{\alpha + \xi _{1}}) \cup A(\frac{1}{\alpha + \xi _{2}})\cup \dots \cup A(\frac{1}{\alpha + \xi _{q^{n}-2}})$ which can be seen as a set of $q^{n}+1$ affine sets. $\langle\sigma^{n}\rangle$ partitions this set of $q^{n}+1$ affine sets. The possible length of orbits are $1$ and $n$. We will consider two cases:  $n=p$ and $n\neq p$.

First suppose that $n=p$. Then $q^n+1\equiv 1 \pmod{n}$ since $q=p^t$. By Theorem \ref{Theorem_Puliya_1}, we know that there is only one affine set fixed under $\langle\sigma^{n}\rangle$ so we conclude that there is one projective linear set fixed under $\langle\sigma^{n}\rangle$ in this case. 

Next suppose that $n\neq p$. We consider two possibilities; $n\mid (q^n-1)$ and $n\nmid (q^n-1)$. Now, if $n\mid (q^n-1)$ then $q^n+1=q^n-1+2\equiv 2 \pmod{n}$. So orbits of length $n$ only are not possible since $n=r>2$. The fact that $q^n+1\equiv 2 \pmod{n}$ implies that a projective linear set fixed under $\langle\sigma^{n}\rangle$ contains $jn+2$ affine sets that are fixed under $\langle\sigma^{n}\rangle$ where $j$ is a non-negative integer. By Theorem \ref{Theorem_Puliya_1}, there are only $n-1$ affine sets fixed under $\langle\sigma^{n}\rangle$. Thus $j =0$. We conclude that there are $\frac{n-1}{2}$ projective linear sets fixed under $\langle\sigma^{n}\rangle$, 2 fixed affine sets in each.

Now suppose that $n\nmid (q^n-1)$. By Theorem \ref{Theorem_Puliya_1} there is no affine set fixed in this case. If $O(\alpha)\in \O$ is fixed under $\langle\sigma^{n}\rangle$ then $\sigma ^{n}(O(\alpha))=O(\alpha)$. So we have $\sigma^{n}(\alpha)=\alpha ^{q^{n}}=[A](\alpha)$ where $A\in GL(2,q^n)$ and $m_A(x)$ is an irreducible quadratic polynomial over $\F_{q^n}$. When we apply $\sigma^{n}$ to $\alpha$ $n$-times we obtain 
$$\alpha=\sigma^{n^2}(\alpha)=[A^n](\alpha)=[I_2](\alpha).$$ We conclude that $A^n=I_2$. Thus $A$ is a matrix of order $n$ over $\mathbb{F}_{q^{n}}$. Next we count the number of such matrices.
We consider two possibilities; $n\nmid (q^n+1)$ and $n\mid (q^n+1)$. 
If $n\nmid (q^n+1)$ matrices of order $n$ do not exist hence there is no projective linear set fixed. Now suppose that $n\mid (q^n+1)$. Then $q^{2n}\equiv 1 \pmod{n}$ so, by Theorem \ref{Lidl} with $d=2$, the factorization of the cyclotomic polynomial $Q_n(x)$ over $\F_{q^n}$ contains $\phi(n)/2=(n-1)/2$ irreducible quadratic factors. Hence, by Theorem \ref{Thm_Yanga}, there are $\frac{q^n(q^n-1)(n-1)}{2}$ matrices of order $n$. For each matrix of order $n$, we have $\alpha^{q^n}=[A](\alpha)$ where $A=\left(
\begin{array}{cc}
a&b\\
c&d
\end{array}\right)\in GL(2,q^n)$ is of order $n$. This gives $\alpha^{q^n}=\frac{a\alpha+b}{c\alpha+d}$. That is, $c\alpha^{q^n+1}+d\alpha^{q^n}-a\alpha-b=0$. So we may assume that $\alpha$ satisfies an equation of the form
\begin{equation}\label{Eqn_Main}
F_n(x)=cx^{q^{n}+1}+dx^{q^n}-ax-b=0,
\end{equation}  
where the coefficients $a,b,c,d \in \F_{q^n}$ come from a matrix of order $n$. Next we count the number of roots of polynomials of the form $F_n(x)$ which lie in $\S$. Let $\S_F$ be the set of roots of all the polynomials $F_n(x)$ which lie in $\S$. Clearly, $|\S_F|$ depends on $a,b,c,d \in \F_{q^n}$ where $ad-bc \neq 0$. We know that the number of matrices in $GL(2,q^n)$ of order $n$ is $\frac{q^n(q^n-1)(n-1)}{2}$. By Theorem \ref{Stich_Bho}, with $s=1$ and $D=n$ we see that $F_n(x)$ factors into polynomials of the same degree $n$. Note that we are taking $s=1$ since $F_n(x)$ is a polynomial over $F_{q^n}$. Thus all $q^n+1$ roots of $F_n(x)$ lie in $\S_F$. As such $|\S_F|=\frac{q^n(q^n-1)(q^n+1)(n-1)}{2}$. Now, since every element of a fixed projective linear set is a root of some polynomial of the form $F_n(x)$ then the number of projective linear sets fixed under $\langle\sigma^{n}\rangle$ in this case is 

$$\frac{|\S_F|}{|O(\alpha)|}=\frac{q^n(q^n-1)(q^n+1)(n-1)}{2q^n(q^n-1)(q^n+1)}=\frac{n-1}{2}.$$

From the foregoing discussion, we have established the following.
\begin{theorem}\label{main_theorem}
%We can define \rho in this way since in this case the orders of the matrices are prime. 
Let $\rho$ be the number of irreducible quadratic polynomials in the factorization of the cyclotomic polynomial $Q_n(x)$ over $\F_{q^n}$. Then the number of projective linear sets fixed under $\langle\sigma^{n}\rangle$ where $n\nmid (q^n-1)$ but $n\mid (q^n+1)$ and all roots of $F_n(x)$ lie in $\S_F$ is $\rho$.
\end{theorem} 
\subsubsection{Action of $\langle\sigma\rangle$, a subgroup of order $n^2$}
Suppose $O(\alpha)\in \O$ is fixed under $\langle\sigma\rangle$. Then $\langle\sigma\rangle$ acts on $O(\alpha)=A(\alpha)\cup A(\frac{1}{\alpha})\cup A(\frac{1}{\alpha + 1})\cup A(\frac{1}{\alpha + \xi _{1}}) \cup A(\frac{1}{\alpha + \xi _{2}})\cup \dots \cup A(\frac{1}{\alpha + \xi _{q^{n}-2}})$ which can be seen as a set of $q^{n}+1$ affine sets. $\langle\sigma\rangle$ partitions this set of $q^{n}+1$ affine sets. The possible length of orbits are $1$, $n$ and $n^2$. We consider the following cases: $n=p$ and $n\neq p$. 

%We know from Section \ref{Action_1_n=r} that there is no affine set fixed under $\langle\sigma\rangle$ if $r\nmid (q-1)$ and $p\neq r$ and that there is one affine set fixed if $p=r$. We also know that there are $r-1$ affine sets fixed if $r\mid (q-1)$.
%We will consider four cases; $n\mid (q^n-1)$, $n\nmid (q^n-1)$ but $p=n$, $n\nmid (q^n-1)$ and $p\neq n$, and $n\nmid (q^n-1)$ but $n\mid (q^n+1)$. 
First suppose that $n=p$. Then $q^n+1\equiv 1 \pmod{n}$ since $q=p^t$. By Theorem \ref{Theorem_Puliya_1}, we know that there is one affine set fixed under $\langle\sigma\rangle$ so we conclude that there is one projective linear set fixed under $\langle\sigma\rangle$. 

Next suppose that $n\neq p$. We consider two possibilities: $n\mid (q^n-1)$ and $n\nmid (q^n-1)$. Now, if $n\mid (q^n-1)$ then $q^n+1=q^n-1+2\equiv 2 \pmod{n}$. So orbits of length $n$ only are not possible since $n=r>2$. As such there must be at least one orbit of length 1, that is, $O(\alpha)$ must contain an affine set that is fixed under $\langle\sigma^{n}\rangle$. By Theorem \ref{Theorem_Puliya_1}, there are $n-1$ affine sets fixed under $\langle\sigma^{n}\rangle$. We conclude that there are $\frac{n-1}{2}$ projective linear sets fixed under $\langle\sigma\rangle$, 2 fixed affine sets in each.

Now suppose that $n\nmid (q^n-1)$.  By Theorem \ref{Theorem_Puliya_1}, we know that there is no affine set fixed under $\langle\sigma\rangle$. So orbits of length 1 are not possible. Additionally, if $n\nmid (q^n+1)$ then it is easy to see that there is no projective linear set fixed. 

If $n\nmid (q^n-1)$ but $n\mid (q^n+1)$ (or $n\mid (q+1)$ since $q^n+1\equiv 0 \pmod{n}$ and $q^n+1\equiv q+1 \pmod{n}$ implies $n\mid (q+1)$) then orbits of length $n$ are possible. By an argument similar to the one in Section \ref{Action_O_n_n=r}, we need to consider how many roots of $F_1(x)=cx^{q+1}+dx^{q}-ax-b=0$ where $a,b,c,d \in \F_{q^n}$ and $ad-bc \neq 0$ lie in $\S$. We see that all $q+1$ roots of $F_1(x)$ lie in $\S_F$ and so by Theorem \ref{main_theorem} the number of projective linear sets fixed under $\langle\sigma\rangle$ given that $n\nmid (q-1)$ but $n\mid (q+1)$ is $\rho$ where $\rho=\frac{n-1}{2}$.  

\subsection{Applying the Cauchy Frobenius Theorem}
Table \ref{Length_2^n+1_n=r} shows the number of projective linear sets which are fixed under the action of various subgroups of $G$. The subgroups are listed in ascending order of the number of elements in the subgroup. We list in column 2 the number of elements in a subgroup which are not already counted in subgroups in the rows above it in the table.
\begin{table}[htbp]
\caption{Number of Projective Linear Sets Fixed}
\label{Length_2^n+1_n=r}
\centering
\begin{tabular}{cccccc}

\hline Subgroup &  No. of elements& No. of fixed       & No. of fixed      & No. of fixed          & No. of fixed\\
        of $G$  & not in previous & projective linear  & projective linear & projective linear     & projective linear\\
                &  subgroup       &sets                & sets              & sets                  & sets\\
                &                 &if $n\neq p$        &if $n\neq p$       &if $n\neq p$           &if $n=p$\\
                &                 &and $n\mid (q^n-1)$ &$n\nmid (q^n-1)$   &$n\nmid (q^n-1)$       &\\
                &                 &                    &and $n\mid (q^n+1)$& and $n\nmid (q^n+1)$  & \\
 
\hline 
\hline $\langle \sigma^{n^2} \rangle$  & 1 & $\frac{q^{n^2}-q^n}{q^{3n}-q^n}$ &$\frac{q^{n^2}-q^n}{q^{3n}-q^n}$&$\frac{q^{n^2}-q^n}{q^{3n}-q^n}$&$\frac{q^{n^2}-q^n}{q^{3n}-q^n}$  \\ 
\hline $\langle \sigma^{n} \rangle$  & $n-1$ & $\frac{n-1}{2}$ & $\frac{n-1}{2}$ & 0& 1\\ 
\hline $\langle \sigma \rangle$  & $n^2-n$ & $\frac{n-1}{2}$ & $\frac{n-1}{2}$& $0$&$1$ \\
\hline   
\end{tabular}  
\end{table}
%\\

By the Cauchy-Frobenius Theorem, we obtain the following result.
\begin{theorem}\label{Thm_all_n=r_q^n+1}
{\rm The number of orbits in $\mathbb{O}$ under the action of $G$ is\\
$$\left\{ 
\begin{array}{r l}
\frac{1}{n^2}\left[\frac{q^{n(n-1)}-1}{q^{2n}-1}+n^2-1\right],& \quad \mbox{if}~n=p \\
\frac{1}{n^2}\left[\frac{q^{n(n-1)}-1}{q^{2n}-1}+\frac{(n^2-1)(n-1)}{2}\right],& \quad \mbox{if}~n\neq p~\mbox{and}~ n\mid q^n-1 \\
\frac{1}{n^2}\left[\frac{q^{n(n-1)}-1}{q^{2n}-1}+\frac{(n+1)(n-1)^2}{2}\right],& \quad \mbox{if}~n\neq p, n\nmid q^n-1~\mbox{and}~n\mid q^n+1\\
\frac{1}{n^2}\left[\frac{q^{n(n-1)}-1}{q^{2n}-1}\right],& \quad \mbox{if}~n\neq p, n\nmid q^n-1 ~\mbox{and}~ n\nmid q^n+1 

\end{array} \right. $$
}
\end{theorem}
\begin{remark}
The number of orbits in $\O$ under the action of $G$ gives us an upper bound on the number of extended irreducible Goppa codes of degree $r$ and length $q^n+1$ where  $n=r$ is prime.

\end{remark}

\begin{example}\label{Ex_2}
If we take $q=2$, $n=r=5$, then there are at most $$\frac{2^{5(5-3)}+2^{5(5-5)}}{25}=\frac{1025}{25}=41$$ irreducible binary Goppa codes of degree $5$ and length $33$.  
\end{example}
The result in Example \ref{Ex_2} was also found in \cite{Magamba_wakale}.

\section{Extended Irreducible Goppa codes of degree $r$ and length $q^n+1$, where $n\neq r$}
Our aim in this section is to obtain the number of extended irreducible Goppa codes of degree $r$ and length $q^n+1$, where $n\neq r$. Before we do that let us first obtain $\S(n,r)$.
\subsection{$\S(n,r)$ where $n$ and $r$ are both prime and $n\neq r$} 
\begin{figure}[htb]
\centering
\begin{tikzpicture}[node distance=1.5cm] 
\node (F1) {$\F_q$}; 
\node (Fn) [above left of=F1]{$\F_{q^{n}}$};
\node (Fr) [above right of=F1]{$\F_{q^{r}}$};
\node (Fnr) [above right of=Fn]{$\F_{q^{nr}}$};

\draw (F1) -- (Fn);
\draw (F1) -- (Fr);
\draw (Fr) -- (Fnr);
\draw (Fn) -- (Fnr);
\end{tikzpicture}
\caption{Lattice of subfields of $\F_{q^{nr}}$}
\label{243}
\end{figure}
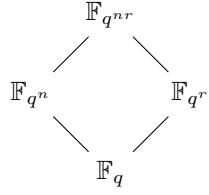

From Figure \ref{243}, elements of $\S(n,r)$ lie in $\F_{q^{nr}}$ and $\F_{q^{r}}$. Hence the number of elements of degree $r$ over $\F_{q^n}$ is $|\S(n,r)|=q^{nr}-q^{n}.$ 

\subsection{Action of $G$ on $\A$}
In this section we find the number of affine sets in $\A$ which are fixed by subgroups of $G$. Thus, we will consider the action of $\langle\sigma\rangle$, $\langle\sigma^{r}\rangle$, $\langle\sigma^{n}\rangle$ and $\langle\sigma^{nr}\rangle$ on $\A$. Clearly, the trivial subgroup $\langle\sigma^{nr}\rangle$ containing the identity fixes every affine set in $\A$. By an argument similar to the one in Section III of \cite{Ryan4} corresponding to the subgroups $\langle\sigma^{n}\rangle$, $\langle\sigma^{r}\rangle$ and $\langle\sigma\rangle$ we obtain the following results.
\begin{theorem}\label{Theorem_Action_n_(n,r)=1}
{\rm The number of affine sets fixed under $\langle\sigma^n\rangle$ when $n\neq r$ is
$$\left\{ 
\begin{array}{r l}
1,& \quad \mbox{if}~ r=p \\
0,& \quad \mbox{if}~ r\neq p~\mbox{and}~r \nmid (q^n-1) \\
r-1,& \quad \mbox{if}~ r\neq p~\mbox{and}~r\mid (q^n-1) 
\end{array} \right. .$$
}
\end{theorem}

\begin{theorem}\label{Action_r_(n,r)=1}
{\rm The number of affine sets fixed by $\langle\sigma^{r}\rangle$ when $n\neq r$ is $$\frac{|\S(1,r)|}{q(q-1)}=\frac{q^{r}-q}{q(q-1)}=\frac{q^{r-1}-1}{q-1}.$$ 
}
\end{theorem}

\begin{theorem}\label{Theorem_Action_1_(n,r)=1}
{\rm The number of affine sets fixed under $\langle\sigma\rangle$ when $n\neq r$ is
$$\left\{ 
\begin{array}{r l}
1,& \quad \mbox{if}~ r=p \\
0,& \quad \mbox{if}~ r\neq p~\mbox{and}~r \nmid (q-1) \\
r-1,& \quad \mbox{if}~ r\neq p~\mbox{and}~r\mid (q-1) 
\end{array} \right. .$$
}
\end{theorem}
\subsection{Action of G on $\O$}
In this section we obtain an upper bound on the number of extended irreducible Goppa codes of degree $r$ and length $q^n+1$ where  $n$ and $r$ are both prime and $n\neq r$. Clearly, the trivial subgroup $\langle\sigma^{nr}\rangle$ containing the identity fixes every projective linear set in $\O$. We now consider the actions of $\langle\sigma^{n}\rangle$, $\langle\sigma^{r}\rangle$ and $\langle\sigma\rangle$. 
\subsubsection{Action of $\langle\sigma^{n}\rangle$, a subgroup of order $r$}\label{ka_sigma_n}
Suppose $O(\alpha)\in \O$ is fixed under $\langle\sigma^{n}\rangle$. Then $\langle\sigma^{n}\rangle$ acts on $O(\alpha)=A(\alpha)\cup A(\frac{1}{\alpha})\cup A(\frac{1}{\alpha + 1})\cup A(\frac{1}{\alpha + \xi _{1}}) \cup A(\frac{1}{\alpha + \xi _{2}})\cup \dots \cup A(\frac{1}{\alpha + \xi _{q^{n}-2}})$ which can be seen as a set of $q^{n}+1$ affine sets. $\langle\sigma^{n}\rangle$ partitions this set of $q^{n}+1$ affine sets. The possible length of orbits are $1$ and $r$. We will consider two possibilities; $r=p$ and $r\neq p$.

First suppose that $r=p$. Then $q^n+1\equiv 1 \pmod{r}$. So orbits of length $r$ only are not possible. So if a projective linear set is fixed under $\langle\sigma^{n}\rangle$ where $r=p$ it has to contain a fixed affine set. Now, by Theorem \ref{Theorem_Action_n_(n,r)=1}, there is one affine set fixed under $\langle\sigma^{n}\rangle$ so it follows that there is one projective linear set fixed.  

Second suppose that $r \neq p$. We consider the following three possibilities 1) $r\mid (q-1)$, 2) $r\mid (q^n-1)$ and $r\nmid (q-1)$ and 3) $r\nmid (q^n-1)$. If $r\mid (q-1)$ and since $q-1\mid(q^n-1)$ for any positive integer $n$ then the fact that $q^n+1=q^n-1+2\equiv 2 \pmod{r}$ implies that orbits of length $r$ only are not possible. It is easy to see that if a projective linear set is fixed under $\langle\sigma^{n}\rangle$ then it contains two fixed affine sets. Now, by Theorem \ref{Theorem_Action_n_(n,r)=1}, there are $r-1$ affine sets fixed under $\langle\sigma^{n}\rangle$. Hence the number of projective linear sets fixed under $\langle\sigma^n\rangle$ given that $r\mid (q-1)$ is $\frac{r-1}{2}$.\\   
Next, suppose that $r\mid (q^n-1)$ but $r\nmid (q-1)$. Then $q^n+1=q^n-1+2\equiv 2 \pmod{r}$. So orbits of length $r$ only are not possible since $r>2$. Now, by Theorem \ref{Theorem_Action_n_(n,r)=1}, there are $r-1$ affine sets fixed under $\langle\sigma^{n}\rangle$. Hence the number of projective linear sets fixed under $\langle\sigma^n\rangle$ given that $r\mid (q^n-1)$ but $r\nmid (q-1)$ is $\frac{r-1}{2}$. \\
Now suppose that $r\nmid (q^n-1)$. Clearly, there is no projective linear set fixed when $r\nmid (q^n-1)$ and $r\nmid (q^n+1)$ since orbits of length $r$ only are not possible and there is no affine set fixed when $r\nmid (q^n-1)$ and $r \neq p$.  However, if $r\nmid (q^n-1)$ but $r\mid (q^n+1)$ then $Q_r(x)$ factors into $\rho=\frac{\phi(r)}{2}$ irreducible quadratic polynomials. By an argument similar to the one in Section \ref{Action_O_n_n=r}, the number of projective linear sets fixed under $\langle\sigma^{n}\rangle$, where there are no affine sets fixed in the decomposition of $O(\alpha)$ and $r\nmid (q^n-1)$ but $r\mid (q^n+1)$ is $\rho=\frac{\phi(r)}{2}=\frac{r-1}{2}$.

\subsubsection{Action of $\langle\sigma^{r}\rangle$, a subgroup of order $n$}\label{Action_of_sigma_r}
Suppose the orbit in $\mathbb{O}$ under the action of $G$ containing $O(\alpha)$ contains $r$ affine sets. Then $O(\alpha)$ is fixed under $\langle\sigma^{r}\rangle$. We claim that each $O(\alpha) \in {\mathbb O}$ fixed under $\langle\sigma^{r}\rangle$ contains precisely $q+1$ affine sets which are fixed under $\langle\sigma^{r}\rangle$. Without loss of generality, suppose $A(\alpha)$ is fixed under $\langle\sigma^{r}\rangle$. Then $A(\alpha)$ contains elements which satisfy the equation $x^{q^{r}}=x$. Assume that $\alpha$ satisfies $x^{q^{r}}=x$. Then it is clear that $\alpha+\nu$ where $\nu \in {\mathbb F}_{q}$ also satisfies the equation $x^{q^r}=x$. So it follows that, for $\frac{\zeta}{\alpha+\nu}+\xi \in A(\frac{1}{\alpha+\nu})$, we have $(\frac{\zeta}{\alpha+\nu}+\xi)^{q^{r}}=\frac{\zeta^{q^{r}}}{\alpha+\nu}+\xi^{q^r} \in A(\frac{1}{\alpha+\nu})$ which implies that $A(\frac{1}{\alpha+\nu})$ is fixed under $\langle\sigma^{r}\rangle$. Since there are $q$ elements in $\mathbb{F}_{q}$ then there are $q$ such affine sets. We now show that no affine set of the form $A(\frac{1}{\alpha+\mu})$, where $\mu\in \mathbb{F}_{q^{n}}\setminus\mathbb{F}_{q}$, in the decomposition of $O(\alpha)$ is fixed. First note that $\mu^{q^{r}}\neq\mu$ for $\mu \in \mathbb{F}_{q^{n}}\setminus\mathbb{F}_{q}$. Since, for $\frac{\zeta}{\alpha+\mu}+\xi \in A(\frac{1}{\alpha+\mu})$ where $\mu \in \mathbb{F}_{q^{n}}\setminus\mathbb{F}_{q}$, we have $(\frac{\zeta}{\alpha+\mu}+\xi)^{q^{r}}=\frac{\zeta^{q^{r}}}{\alpha+\mu^{q^{r}}}+\xi^{q^{r}}\notin A(\frac{1}{\alpha+\mu})$ then $A(\frac{1}{\alpha+\mu})$ is not fixed under $\langle\sigma^{r}\rangle$. So we conclude that there are $q+1$ affine sets in $O(\alpha)$ which are fixed under $\langle\sigma^{r}\rangle$. By Section \ref{Action_r_(n,r)=1}, there are $\frac{q^{r-1}-1}{q-1}$  affine sets fixed under $\langle\sigma^{r}\rangle$. Hence the number of projective linear sets in ${\mathbb O}$ which are fixed under $\langle\sigma^{r}\rangle$ is $$\frac{q^{r-1}-1}{(q-1)(q+1)}=\frac{q^{r-1}-1}{q^2-1}.$$

\subsubsection{Action of $\langle\sigma\rangle$, a subgroup of order $nr$}
Suppose $O(\alpha)\in \O$ is fixed under $\langle\sigma\rangle$. Then $\langle\sigma\rangle$ acts on $O(\alpha)=A(\alpha)\cup A(\frac{1}{\alpha})\cup A(\frac{1}{\alpha + 1})\cup A(\frac{1}{\alpha + \xi _{1}}) \cup A(\frac{1}{\alpha + \xi _{2}})\cup \dots \cup A(\frac{1}{\alpha + \xi _{q^{n}-2}})$ which can be seen as a set of $q^{n}+1$ affine sets. $\langle\sigma\rangle$ partitions this set of $q^{n}+1$ affine sets. The possible length of orbits are $1$, $n$, $r$ and $nr$. We will consider two possibilities: $r=p$ and $r\neq p$. 

%\subsubsection*{Projective linear sets fixed under $\langle\sigma\rangle$ when $p=r$}
Suppose that $r=p$. Then $q^n+1\equiv 1 \pmod{r}$. So orbits of length $r$ only are not possible. So if a projective linear set is fixed under $\langle\sigma\rangle$ it has to contain a fixed affine set. Now, by Theorem \ref{Theorem_Action_1_(n,r)=1}, there is only one affine set fixed under $\langle\sigma\rangle$ so it follows that there is only one projective linear set fixed under $\langle\sigma\rangle$. Now if $O(\alpha)$ fixed under $\langle\sigma\rangle$ contains an orbit of length $n$ then $O(\alpha)$ is also fixed under $\langle\sigma^n\rangle$. But we know that if $r=p$ a projective linear set fixed under $\langle\sigma^n\rangle$ contains one fixed affine set so there is no projective linear set fixed in this case. Next we consider the possibility of a fixed projective linear set where there are orbits of length $n$, $r$ and $nr$. That is, we can find non-negative integers $x$ and $y$ and $z$ such that $nx+ry+nrz=q^n+1$. Since such a projective linear set contains an orbit of length $n$, it is also fixed under $\langle\sigma^n\rangle$. By Section \ref{ka_sigma_n}, a projective linear set fixed under $\langle\sigma^n\rangle$ where $r=p$ contains one fixed affine set. We conclude that if a projective linear set is fixed under $\langle\sigma\rangle$ where $r=p$ then it contains one affine set fixed. Hence there is one projective linear set fixed under $\langle\sigma\rangle$. 

%\subsubsection*{Projective linear sets fixed under $\langle\sigma\rangle$ when $p\neq r$}
Now suppose that $r \neq p$. We consider the following three possibilities 1) $r\mid (q-1)$, 2) $r\mid (q^n-1)$ but $r\nmid (q-1)$ and 3) $r\nmid (q^n-1)$. 

We begin by looking at the case $r\mid (q-1)$. For any positive integer $n$, $q-1\mid (q^n-1)$ so $q^n+1=q^n-1+2\equiv 2 \pmod{r}$. As such, orbits of length $r$ only are not possible since $r>2$. The fact that $q^n+1\equiv 2 \pmod{r}$ implies that a projective linear set fixed under $\langle\sigma\rangle$ contains $jr+2$ affine sets that are fixed under $\langle\sigma\rangle$ where $j$ is a positive integer. Now, by Theorem \ref{Theorem_Action_1_(n,r)=1}, there are only $r-1$ affine sets fixed under $\langle\sigma^{n}\rangle$. This implies that $j=0$. So there is no projective linear set fixed. However, if a projective linear set is fixed under $\langle\sigma\rangle$ and contains an orbit of length $n$, then $\langle\sigma^n\rangle$ fixes such a projective linear set. By Section \ref{ka_sigma_n}, a projective linear set fixed under $\langle\sigma^n\rangle$ contains two fixed affine sets. Hence the number of projective linear sets fixed under $\langle\sigma\rangle$ given that $r\mid (q-1)$ is $\frac{r-1}{2}$. 

%The fact that $q^n+1=q^n-1+2\equiv 2 \pmod{n}$ implies that if a projective linear set is fixed under $\langle\sigma\rangle$ then it contains 2 affine sets that are fixed under $\langle\sigma\rangle$. 

Now, suppose that $r\mid (q^n-1)$ where $r\nmid (q-1)$. By Theorem \ref{Theorem_Action_1_(n,r)=1}, there are no affine sets fixed, so the possible lengths of an orbit are $n$, $r$ and $nr$. We consider two cases 1) $r\nmid (q+1)$ and 2) $r\mid (q+1)$.
Suppose $r\nmid (q+1)$, then orbits of length $r$ only are not possible. So there is no projective linear set fixed. Next we consider orbits of length $n$. We know that $q^n+1\equiv q+1\pmod{n}$, by Fermat's Little Theorem. So if $n\mid (q+1)$ then orbits of length $n$ are possible. Now, a projective linear set fixed under $\langle\sigma\rangle$ containing an orbit of length $n$ is also fixed under $\langle\sigma^n\rangle$. From Section \ref{ka_sigma_n}, we know that a projective linear set fixed under    $\langle\sigma^n\rangle$ contains 2 fixed affine sets and $q^n-1$ affine sets that are permuted in orbits of length $r$. Now, since $n\neq r$ there is no projective linear set fixed in this case. 
%Suppose that $r\mid (q^n-1)$ where $r\nmid (q-1)$ and $r\nmid (q+1)$. If $O(\alpha)$ is fixed under $\langle\sigma\rangle$  then $\sigma(\alpha)=\alpha^q=\frac{a\alpha+b}{c\alpha+b}$. Now, applying $\sigma$ to $\alpha$ $r$-times we obtain $\alpha^{q^r}=\alpha=[A^r](\alpha)=[I_2](\alpha)$. This implies that $A\in GL(2,q)$ is of order $r$. Observe that the fact that $r\nmid (q-1)$ and $r\nmid (q+1)$ implies that matrices of order $r$ over $\F_q$ do not exist and as such no projective linear set can be fixed under $\langle\sigma\rangle$ given these parameters.    

%Now if $\langle\sigma^{n}\rangle$ fixes $O(\alpha)$ then $\alpha^{q^n}=\frac{a\alpha +b }{c\alpha +d}$ where $a,b,c,d \in \F_{q^n}$ and $ad-bc \neq 0$, see Section \ref{Actions}. $\alpha^{q^n}=\frac{a\alpha +b }{c\alpha +d}$ implies that $c\alpha^{q^n+1}+d\alpha^{q^n}=a\alpha +b$. Thus, we may assume that $\alpha$ satisfies an equation of the form
%
%\begin{equation}\label{Eqn_Main}
%F_n(x)=cx^{q^{n}+1}+dx^{q^n}-ax-b=0.
%\end{equation}
%
%
%  Clearly orbits of length $r$ only do not exist since $q^n+1=q^n-1+2\equiv 2 \pmod{r}$ implies that $r\nmid (q^n+1)$. Also, by an argument similar to the one in the paragraph above, an $O(\alpha)$ fixed under $\langle\sigma\rangle$ that contains orbits of length $n$ only are not possible. From this discussion we see that it is not possible to find a projective linear set fixed under $\langle\sigma\rangle$ where there are orbits of length $n$, $r$ and $nr$ such that for positive integers $x$ and $y$ and an integer $z\geq 0$, $nx+ry+nrz=q^n+1$.

However, if $r\mid (q^n-1)$ where $r\nmid (q-1)$ and $r\mid (q+1)$ then matrices $A\in GL(2,q)$ of order $r$ exist. By an argument similar to the one in Section \ref{Action_O_n_n=r}, we need to find how many roots of $F_1(x)=cx^{q+1}+dx^{q}-ax-b=0$ where $a,b,c,d \in \F_{q}$ and $ad-bc \neq 0$ lie in $\S$. We conclude that the number of projective linear sets fixed under $\langle\sigma\rangle$ is $\frac{r-1}{2}$. 
%Note that if $r\mid (q^n-1)$ where $r\nmid (q-1)$ and $r\nmid (q+1)$ then there are no projective linear sets fixed under $\langle\sigma\rangle$.   

Next, suppose that $r\nmid (q^n-1)$. By Theorem \ref{Theorem_Action_1_(n,r)=1} there is no affine set fixed under $\langle\sigma\rangle$. As such, the possible lengths of an orbit are $n$, $r$ and $nr$. If $r\nmid (q^n-1)$ and $r\nmid (q^n+1)$ then matrices $A\in GL(2,q^n)$ of order $r$ do not exist. Similarly, if $r\nmid (q^n-1)$ and $r\nmid (q+1)$ then matrices $A\in GL(2,q)$ of order $r$ do not exist. In either case, there are no projective linear sets fixed under $\langle\sigma\rangle$. However, if $r\nmid (q^n-1)$ and $r\mid (q+1)$ then, as above, we need to find how many roots of $F_1(x)=cx^{q+1}+dx^{q}-ax-b=0$, where $a,b,c,d \in \F_{q}$ and $ad-bc \neq 0$, lie in $\S$. By an argument similar to the one in Section \ref{Action_O_n_n=r}, the number of projective linear sets fixed under $\langle\sigma\rangle$ is $\frac{r-1}{2}$. 

Lastly, suppose that there is a possibility of having a projective linear set $O(\alpha)$  fixed  under $\langle\sigma\rangle$ with a combination of different orbit lengths. That is, we can find non-negative integers $x,y$ and $z$ all not equal to zero such that $nx+ry+nrz=q^n+1$. Observe that if $y\neq 0$ then such a projective linear set is also fixed under $\langle\sigma^r\rangle$. Now, a projective linear set fixed under $\langle\sigma^r\rangle$ contains $q+1$ fixed affine sets, see Section \ref{Action_of_sigma_r}. This implies $r=q+1$. We see that we have already dealt with this case. Moreover, if $x\neq 0$ then $O(\alpha)$ is also fixed under $\langle\sigma^n\rangle$. By Section \ref{ka_sigma_n}, if there is no fixed affine set in the decomposition of $O(\alpha)$ then $r\mid (q^n+1)$ and $\langle\sigma^n\rangle$ fixes $\frac{r-1}{2}$ projective linear sets. This result is consistent with the results above.       

\subsection{Applying the Cauchy-Frobenius Theorem}
Tables \ref{Length_q^n+1_(n,r)=1_C}, \ref{Length_q^n+1_(n,r)=1_D} and \ref{Length_q^n+1_(n,r)=1_D_1} show the number of projective linear sets which are fixed under the action of various subgroups of $G$.
\begin{table}[htbp]
\caption{Number of Projective Linear (PL) Sets fixed}
\label{Length_q^n+1_(n,r)=1_C}
\centering
\begin{tabular}{ccccccc}
\hline Subgroup& No. of elements & No. of fixed & No. of fixed  & No. of fixed      & No. of fixed & No. of fixed\\
        of $G$ & not in previous & PL sets      & PL sets       & PL sets           & PL sets& PL sets\\
               & subgroup        & if $r=p$     & if $r\neq p$  & if $r\neq p$,& if $r\neq p$,& if $r\neq p$,\\
               &                 &              & and $r\mid (q-1)$ &$r\mid (q^n-1)$     & $r\nmid (q^n-1)$,&  $r\nmid (q^n-1)$\\
               &                 &              &               & and $r\nmid (q-1)$  & and $r\mid (q^n+1)$ & and $r\nmid (q^n+1)$ \\
\hline 
\hline $\langle \sigma^{nr} \rangle$ &  1 & $\frac{q^{nr}-q^n}{q^{3n}-q^n}$& $\frac{q^{nr}-q^n}{q^{3n}-q^n}$&$\frac{q^{nr}-q^n}{q^{3n}-q^n}$&$\frac{q^{nr}-q^n}{q^{3n}-q^n}$&$\frac{q^{nr}-q^n}{q^{3n}-q^n}$ \\ 
\hline $\langle \sigma^{n} \rangle$  & $r-1$ & 1&$\frac{r-1}{2}$&$\frac{r-1}{2}$&$\frac{r-1}{2}$&0  \\ 
\hline $\langle \sigma^{r} \rangle$  & $n-1$ & $\frac{q^{r-1}-1}{q^2-1}$& $\frac{q^{r-1}-1}{q^2-1}$& $\frac{q^{r-1}-1}{q^2-1}$&$\frac{q^{r-1}-1}{q^2-1}$ &$\frac{q^{r-1}-1}{q^2-1}$\\ 
%\hline $\langle \sigma \rangle$  & $(n-1)(r-1)$ & 1&$\frac{r-1}{2}$&$\rho^{\prime}$&$\rho^{\prime\prime}$&0\\
\hline 

\end{tabular}  
\end{table}

\begin{table}[htbp]
\caption{Number of Projective Linear (PL) Sets fixed}
\label{Length_q^n+1_(n,r)=1_D}
\centering
\begin{tabular}{ccccc}
\hline
 Subgroup& No. of elements & No. of fixed & No. of fixed      & No. of fixed \\
        of $G$ & not in previous & PL sets      & PL sets           & PL sets  \\
               & subgroup        & if $r=p$     & if $r\neq p$      & if $r\neq p$,\\
               &                 &              & and $r\mid (q-1)$ &$r\mid (q^n-1)$,\\
               &                 &              &                   & $r\nmid (q-1)$,\\
               &                 &              &                   & and $r\mid (q+1)$ \\
\hline 

\hline $\langle \sigma \rangle$  & $(n-1)(r-1)$ & 1&$\frac{r-1}{2}$&$\frac{r-1}{2}$\\
\hline 

\end{tabular}  
\end{table}

\begin{table}[htbp]
\caption{Number of Projective Linear (PL) Sets fixed}
\label{Length_q^n+1_(n,r)=1_D_1}
\centering
\begin{tabular}{ccccc}
\hline Subgroup& No. of elements &No. of fixed      & No. of fixed       & No. of fixed\\
        of $G$ & not in previous &PL sets         & PL sets            & PL sets\\
               & subgroup        &if $r\neq p$,     & if $r\neq p$,      & if $r\neq p$,\\
               &                 &$r\nmid (q^n-1)$, & $r\nmid (q^n-1)$,  & $r\nmid (q^n-1)$,\\
               &                 &$r\mid (q^n+1)$,    & and $r\mid (q+1)$  & and $r\nmid (q^n+1)$\\
               &                 & and $r\nmid (q+1)$  &                    &  \\
\hline 

\hline $\langle \sigma \rangle$  & $(n-1)(r-1)$ &0&$\frac{r-1}{2}$&0\\
\hline 

\end{tabular}  
\end{table}

\begin{theorem}\label{Thm_all_(n,r)=1_q^n+1}
The number of orbits in $\mathbb{O}$ under the action of $G$ is 
%\begin{enumerate}
%\item $\begin{array}{cc}
%\frac{1}{nr}\left[\frac{q^{n(r-1)}-1}{q^{2n}-1}+n(r-1)+\frac{(n-1)(q^{r-1}-1)}{q^2-1}\right],&\quad \mbox{if}~ r=p.
%\end{array}$
%\item $\begin{array}{cc}\frac{1}{nr}\left[\frac{q^{n(r-1)}-1}{q^{2n}-1}+\frac{n(r-1)^2}{2}+\frac{(n-1)(q^{r-1}-1)}{q^2-1}\right],\quad ~\mbox{if}~ r\neq p, ~r\mid (q-1)~\mbox{or}~r\neq p, ~r\mid (q^n-1), ~r\nmid (q-1),~ r\mid (q+1) ~\mbox{or}~ r\neq p,~ r\nmid (q^n-1),~ r\mid (q^n+1)~\mbox{and}~ r\mid (q+1).\end{array}$
%
%\item $\frac{1}{nr}\left[\frac{q^{n(r-1)}-1}{q^{2n}-1}+\frac{(r-1)^2}{2}+\frac{(n-1)(q^{r-1}-1)}{q^2-1}\right],\quad\quad \mbox{if}~ r\neq p,~ r\nmid (q^n-1),r\mid (q^n+1), ~r\nmid (q+1).$
%
%
%\item $\frac{1}{nr}\left[\frac{q^{n(r-1)}-1}{q^{2n}-1}+\frac{(n-1)(q^{r-1}-1)}{q^2-1}\right],\quad\quad\quad\quad\quad\quad \mbox{if}~ r\neq p,~ r\nmid (q^n-1),r\nmid (q^n+1).$
%
%\end{enumerate}
\begin{enumerate}
\item $\begin{array}{rl}
\frac{1}{nr}\left[\frac{q^{n(r-1)}-1}{q^{2n}-1}+n(r-1)+\frac{(n-1)(q^{r-1}-1)}{q^2-1}\right],&\mbox{if}~ r=p.
\end{array}$
\item $\begin{array}{rl}
\frac{1}{nr}\left[\frac{q^{n(r-1)}-1}{q^{2n}-1}+\frac{n(r-1)^2}{2}+\frac{(n-1)(q^{r-1}-1)}{q^2-1}\right],&\mbox{if}~ r\neq p~\mbox{and}~r\mid (q-1)~\mbox{or}~r\neq p, ~r\mid (q^n-1), ~r\nmid (q-1) ~\mbox{and}~r\mid (q+1)\\
&\mbox{or}~ r\neq p,~ r\nmid (q^n-1),~ r\mid (q^n+1)~\mbox{and}~ r\mid (q+1).\end{array}$

\item $\begin{array}{rl}
\frac{1}{nr}\left[\frac{q^{n(r-1)}-1}{q^{2n}-1}+\frac{(r-1)^2}{2}+\frac{(n-1)(q^{r-1}-1)}{q^2-1}\right],&\mbox{if}~ r\neq p,~ r\nmid (q^n-1),r\mid (q^n+1)~\mbox{and}~r\nmid (q+1).
\end{array}
$

\item $\begin{array}{rl}
\frac{1}{nr}\left[\frac{q^{n(r-1)}-1}{q^{2n}-1}+\frac{(n-1)(q^{r-1}-1)}{q^2-1}\right],&\mbox{if}~ r\neq p,~ r\nmid (q^n-1)~\mbox{and}~r\nmid (q^n+1).
\end{array}$

\end{enumerate}

\end{theorem}

\begin{remark}
The number of orbits in $\O$ under the action of $G$ gives us an upper bound on the number of extended irreducible Goppa codes of degree $r$ and length $q^n+1$ where $q=p^t$, $n$ and $r$ are prime numbers and $n\neq r$.
\end{remark}

\begin{example}\label{Ex_4}
If we take $q=2$, $n=11$ and $r=5$, then there are at most $76,261$ extended irreducible binary Goppa codes of degree $5$ and length $2,049$.  
\end{example}

\section{Conclusion}
In this paper we have produced an upper bound on the number of extended irreducible Goppa codes of degree $r$ and length $q^n+1$ where $n$ and $r$ are both prime numbers.

\end{document}